# RAMANUJAN'S THOUGHTS FROM GOD

Frank Aiello

**ABSTRACT.** We review some of Ramanujan's contributions to mathematics. These include his $1/\pi$ series, his work on modular forms, and his work on partitions. We very briefly review his life, including his collaboration with Hardy. Finally, we discuss what any prospective mathematician should know about Ramanujan's work, including the rich relationship between his work on partitions and his work on modular forms. The title is a reference to a direct quote from Ramanujan himself: an equation means nothing to me unless it expresses a thought of God.

## Introduction

Ramanujan's work was awe-inspiring in terms of its sheer creativity. What distinguishes Ramanujan's work from those of his predecessors is that Ramanujan wrote many theorems, but he proved very few. Ramanujan was someone whose mathematical insights were inspired by intuition. I would argue that this style of doing mathematics was largely a result of the fact that Ramanujan taught himself most of the mathematics he knew through G.S. Carr's *A Synopsis of Elementary Results in Pure Mathematics*. This book contained many theorems but was short on proofs, and was in fact outdated by the time Ramanujan came to use it.

## His Work

Ramanujan came up with many important breakthroughs and these were all inspired by intuition in one way or another. Some of these breakthroughs include: Ramanujan's series for $1/\pi$, his work on mock theta functions, and Ramanujan's congruences for the partition function p(n). We start with one of Ramanujan's most remarkable discoveries: his series for $1/\pi$. Ramanujan came up with the following formula for the approximation of pi:

$$\frac{1}{\pi} = \sum_{n=0}^{\infty} \frac{\sqrt{8}(4n)!\,(1103 + 26390n)}{9801(n!)^4 396^{4n}}$$

The more interesting question is Ramanujan's motivation for such a formula. Actually, the formula is a special case for N = 58 of the following theorem:

$$\frac{1}{\pi} = \sum_{n=0}^{\infty} \frac{\left(\frac{1}{4}\right)_n \left(\frac{1}{2}\right)_n \left(\frac{3}{4}\right)_n d_n(N)}{(n!)^3} x_N^{2n-1}$$

where,

$$x_N := \frac{4k_N(k'_N)^2}{(1+k^2)} := \left(\frac{g_N^{12} + g_N^{-12}}{2}\right)^1$$

with

$$d_n(N) = \left[\frac{\alpha(N)x_N^{-1}}{1+k^2{}_N} - n\sqrt{n}\left(\frac{g_N^{12} + g_N^{-12}}{2}\right)^1\right]$$

and

$$k_N := k(e^{-\pi\sqrt{N}}), \quad g_N^{12} = \frac{(k'_N)^2}{(2k_N)}$$

This formula was originally written down in Ramanujan's paper "Modular Equations and Approximations to π". The formula given is a 58th order modular equation. In fact, the derivation of this formula is based on the theory of complete elliptic integrals. The complete integral of the first kind is defined as:

$$K(k) = \int_0^{\frac{\pi}{2}} (1 - k^2 \sin^2\theta)\, d\theta$$

The complete integral of the second kind is defined as:

$$E(k) = \int_0^{\frac{\pi}{2}} (1 - k^2 \sin^2\theta)\, d\theta$$

In fact, there is a close connection between how elliptic integrals transform and the very rapid approximation of pi. This connection was first elucidated by Ramanujan in his 1914 paper "Modular Equations and Approximations to π". Ramanujan came up with many 1/π series and these series could be used to approximate pi remarkably well, such as being able to compute pi to a billion decimal places.

Besides his $1/\pi$ series, another of Ramanujan's great contributions to mathematics was his work on mock theta functions. These mock theta functions were defined on Ramanujan's deathbed. Mock theta functions are defined as incomplete harmonic Maas forms. Recall Euler's generating function:

$$\sum_{n=0}^{\infty} p(n)q^n = \prod_{n=1}^{\infty} \frac{1}{1-q^n} = 1 + q + 2q^2 + 3q^3 + 5q^4 + 7q^5 + \cdots$$

This is what mathematicians refer to as a weakly holomorphic modular form. A holomorphic modular form f of integer weight k, level N, and character $\chi$, over a group G, is a complex function defined on the upper half-plane H which satisfies the transformation $f(\gamma z) = \chi(d)(cz + d)^k f(z)$ for any $z \in H$ and $\gamma \in G$, where G is any subgroup of $SL_2(Z)$ containing the principal congruence subgroup $\Gamma(N)$, where f must be holomorphic on H and at all cusps of G.

A weakly holomorphic modular form is similar to a holomorphic modular form, except that it is allowed to have poles at cusps. Ramanujan's work on mock theta functions fits properly in the theory of modular forms, where a modular form of weight k and level N is an analytic function f from the complex upper half plane to the complex plane such that

$$f\left(\frac{az+b}{cz+d}\right) = (cz+d) * f(z)$$

For all integers a, b, c, d with $ad - bc = 1$

Finally, let us touch on Ramanujan's work on the partition function. Much of Ramanujan's work in the theory of modular forms (and mock theta functions in particular) was a result of his work on the asymptotic properties of the partition function p(n). A partition of a positive integer n is a unique way of writing n as a sum of integers. For example, there are 15 partitions of 7: 7, 6+1,

5+2, 5+1+1, 4+3, 4+2+1, 4+1+1+1, 3+3+1, 3+2+2, 3+2+1+1, 3+1+1+1+1, 2+2+2+1, 2+2+1+1+1, 2+1+1+1+1+1, 1+1+1+1+1+1+1.

Ramanujan collaborated with Hardy on the *Hardy-Ramanujan asymptotic partition formula*. For n a positive integer, let p(n) denote the number of unordered partitions of n; then the value of p(n) is given asymptotically by:

$$p(n) \sim \frac{1}{4n\sqrt{3}} e^{\tau\sqrt{n/6}}$$

Besides his work on the asymptotic partition formula, Ramanujan came up with three congruences for the partition function p(n). These are:

p(5n + 4) ≡ 0 (mod 5)

p(7n + 5) ≡ 0 (mod 7)

p(11n + 6) ≡ 0 (mod 11)

where n is a nonnegative integer.

The proofs of these three congruences for the partition function p(n) rests on the Ramanujan-Eisenstein series:

$$P(q) := 1 - 24 \sum_{k=1}^{\infty} \frac{kq^k}{1-q^k}, |q| < 1$$

$$Q(q) := 1 + 240 \sum_{k=1}^{\infty} \frac{k^3 q^k}{1-q^k}, |q| < 1$$

$$R(q) := 1 - 504 \sum_{k=1}^{\infty} \frac{k^5 q^k}{1-q^k}, |q| < 1$$

In fact, there was one more congruence Ramanujan came up with: $p(13n - 7) \equiv 11\tau(n) \pmod{13}$. We have very briefly reviewed three of Ramanujan's major contributions to mathematics: his $1/\pi$ series, his work on mock theta functions, and his work on partitions. The first of these results is based on the theory of elliptic functions, the second fits properly into his theory of modular forms, and the final result – the proof of his congruences for the partition function $p(n)$ - is based on the Ramanujan-Eisenstein series.

## His Life

Although Ramanujan could easily be compared to Newton, Jacobi, Euler, Gauss and many of the other great mathematicians throughout history, what arguably set Ramanujan apart from these mathematicians was his life and his way of thinking about mathematics. He was almost entirely self-taught, and arguably it was his remarkable intuition for mathematics that allowed him to understand advanced mathematics with almost no formal training. He entered primary school in Kumbakonam, although he attended several different primary schools before attending the Town High School in Kumbakonam in 1898.

He entered Pachiyappa's College in 1906, one of the best colleges in Southern India. His intention was to pass the Fine Arts examination which would allow him to be admitted to the University of Madras. He attended lectures at Pachiyappa's College but became ill after three months of study. He took the Fine Arts examination but failed every subject but mathematics, thus failing the exam. This obviously meant that he could not enter the University of Madras.

Ramanujan continued his mathematical work, including studying continued fractions in 1908. Ramanujan discovered about 200 continued fractions. One of these continued fractions is:

$$\cfrac{1}{1+\cfrac{e^{-2\pi}}{1+\cfrac{e^{-4\pi}}{1+\cfrac{e^{-6\pi}}{1+\cdots}}}} = \left(\sqrt{\frac{5+\sqrt{5}}{2}} - \frac{\sqrt{5}+1}{2}\right)e^{2\pi/5} = 0.9981360\ldots$$

Other of Ramanujan's continued fractions include several continued fractions for pi such as:

$$\pi = \cfrac{4}{1+\cfrac{1^2}{2+\cfrac{3^2}{2+\cfrac{5^2}{2+\cdots}}}}$$

and

$$\pi = 3 + \cfrac{1^2}{6+\cfrac{3^2}{6+\cfrac{5^2}{6+\cdots}}}$$

In 1910, Ramanujan tutored various students in mathematics in the city of Madras and he also walked around the city offering to do accounting work for various businesses. In 1911, Ramanujan published a paper on the Bernoulli numbers; he discovered that the denominators of the fractions of Bernoulli numbers were always divisible by 6. Ramanujan's paper on Bernoulli numbers was published in the *Journal of the Indian Mathematical Society*. Once again, Ramanujan's work on the Bernoulli numbers, like his work on continued fractions, was ultimately driven by his strong mathematical intuition.

Ramanujan obtained a clerical job in 1912. Ramaswamy Aiyer was the Chief Accountant at the Madras Port Trust Office. Ramanujan attempted to find a job at the government revenue department; Ramanujan did not have a resume to show Aiyer but he spoke of his mathematical interests and showed Aiyer a notebook containing his mathematical results.

Aiyer had no mind to smother Ramanujan's genius by an appointment in the lowest rungs of the revenue department. Instead, Aiyer contacted the secretary of the Indian Mathematical Society, R. Ramachandra Rao. Despite initial resistance, after meeting with Ramanujan in person, he was convinced that he was dealing with a mathematical genius. Rao agreed to provide financial support for Ramanujan, and Aiyer began publishing his work in the *Journal of the Indian Mathematical Society*. In March 1912 Ramanujan secured a position as an accounting clerk with the Madras Port Trust.

Ramanujan wrote to G.H. Hardy in January of 1913. Hardy recognized the importance of Ramanujan's work, but demanded that he see proofs of Ramanujan's theorems, which were sorely lacking. Despite the fact that Ramanujan produced remarkable results on elliptic integrals, hypergeometric series, and functional equations of the Riemann zeta function, he had only the vaguest notion of mathematical proof.

Hardy spent two hours going over Ramanujan's work with another mathematician, J.E. Littlewood. Hardy and Littlewood were at a roadblock because they couldn't tell whether they were dealing with the work of a crackpot or a genius. They eventually concluded that it must be a work of genius. To quote Hardy in referring to Ramanujan's theorems: "a single look at them is enough to show that they could only be written by a mathematician of the highest class. They must be true, because if they were not true, no one would have had the imagination to invent them."

Ramanujan traveled to Cambridge in 1914. There, he began collaborating with Hardy and Littlewood. Hardy and Ramanujan eventually published a groundbreaking paper in 1917 in which they applied the Farey dissection of the unit circle to obtain their asymptotic formula for the number of partitions of an integer. Over the next few decades, Hardy and Littlewood developed this circle method into a powerful tool for tackling problems in additive number theory.

As we have already seen, Ramanujan is someone who was the master of mathematical intuition. There were times, however, where this intuition led him astray. An example was Ramanujan's formulas for the distribution of primes. Ramanujan had thought he had found an exact formula for the prime counting function $\pi(n)$, the number of primes not exceeding n. Upon closer inspection by Hardy and Littlewood, this formula was shown to be incorrect.

To quote Hardy: "Ramanujan's theory of primes was vitiated by his ignorance of the theory of functions of a complex variable. It was (so to say) what the theory might be if the zeta function had no complex zeroes. His method depended on the wholesale use of divergent series…that his proofs should be invalid was only to be expected. But the mistakes went deeper than that, and many of the actual results were false. He had obtained the dominant terms of the classical formulae, although by invalid methods; but none of them are such close approximations as he supposed."

In Cambridge, Ramanujan continued to work on partitions, elliptic functions, and continued fractions. Ramanujan did his thesis on highly composite numbers and was awarded a B.A. in March 1916. In 1916, Ramanujan's health was declining again. This was largely his own fault as he would refuse to eat or sleep, working on mathematics obsessively. In 1916, he released a groundbreaking paper titled "On certain arithmetical functions". In this paper, Ramanujan investigated the properties of the Fourier coefficients of modular forms. Also in this paper

Ramanujan enunciates three fundamental conjectures that served as the guiding force for the development of his theory of modular forms.

Ramanujan developed the notion of a *mock modular form*. As we've already stated, mock theta functions are reminiscent of modular forms, which are, loosely speaking, holomorphic functions on the upper half plane equipped with certain symmetries. Ramanujan actually called this notion a "mock theta function" because he thought of his functions as pseudo-modular forms. As an example, one of Ramanujan's mock theta functions is given by:

$$f(q) := \sum n \geq 0 \, q^{n^2} / (-q; q)_n^2$$

where

$$(a; q)_n := \prod_{j=0}^{n-1}(1 - aq^j).$$

These mock theta functions have applications in physics, where they are used to determine critical dimensions in some string theories. By far one of the most important results to come out of Ramanujan's work on modular forms was Ramanujan's conjecture. Ramanujan's conjecture concerns the estimation of Fourier coefficients of the weight 12 holomorphic cusp form Δ for SL(2, Z) on the upper half plane H. It states:

$$\Delta(z) = (2\pi)^{12} e^{2\pi i z} \prod_{n=1}^{\infty}(1 - e^{2\pi i n z})^{24} = (2\pi)^{12} \sum_{n=1}^{\infty} \tau(n) e^{2\pi i n z}$$

where z ∈ H = {z ∈ C| Im z > 0}. By 1917, Ramanujan's health was deteriorating fast. He fell into a deep depression and attempted suicide by attempting to jump on a set of train tracks. Ramanujan was hospitalized in London and was said to have tuberculosis but it is more likely this was to cover up the failed suicide attempt. Hardy took a cab to visit Ramanujan in the

hospital; the only thing he could think to mention was that the number of his taxicab was 1729, a dull number. On the contrary, Ramanujan said it is quite an interesting number. It is the smallest number expressible as a sum of two cubes in two distinct ways:

$$1729 = 12^3 + 1^3 = 10^3 + 9^3$$

Ramanujan was in sanitoria at Wells, Matlock, and in London and it was not until autumn of 1918 that he showed any improvement. In early 1919 it seemed he had recovered sufficiently for the voyage home to India. In 1919, Ramanujan published a new proof of Bertrand's postulate. This postulate states for every $n > 1$, there is always at least one prime $p$ such that $n < p < 2n$. From this paper also came the notion of the Ramanujan prime. The nth Ramanujan prime is the smallest prime such that there are at least n primes between x and 2x for any x such that $2x > p$. Thus, the first few Ramanujan primes are 2, 11, 17, and so on. Unfortunately, Ramanujan's improvement was short-lived. At the end of March 1919, he had a relapse and his health deteriorated tirelessly until his death in April 1920.

Despite having no formal training in mathematics, Ramanujan's mathematical insights stand out for both their profundity and creativity. Truly, it was very difficult for any mathematician to teach Ramanujan any formal mathematics for the simple reason that any attempt to do so led to a flood of new ideas, each more original than the last. In the history of mathematics, Ramanujan was and will forever be remembered as the master of mathematical intuition. Hardy described his collaboration with Ramanujan as the one romantic incident of his life.

## What Prospective Mathematicians Should Know

Ramanujan's mathematical works have had much influence on mathematics education. As we noted, much of Ramanujan's work is based on Eisenstein series. In graduate work, one typically first becomes acquainted with Eisenstein series from their association with the modular group. The modular group G is the subgroup SL(2, Z)/{±1} in PSL(2, R), consisting of matrices with coefficients in Z up to equivalence by ±1. Eisenstein series for the modular group come in two types. One is a series in which the summation is over cosets of a group, and the other is a series in which the summation is taken over elements of a lattice. These two types of series essentially coincide.

Let $\chi$ and $\psi$ be Dirichlet characters mod L and mod M, respectively. For any positive integer $k \geq 3$, we put:

$$(1)\ E_k(z; \chi, \psi) = \sum_{m,n=-\infty}^{\infty} \chi(m)\psi(n)(mz+n)^{-k}$$

Where $\Sigma$ is the summation over all pairs of integers (m, n) except (0, 0). It can be shown that there is a real-valued continuous function r(z) satisfying:

$$(2)\ \sum_{m,n=-\infty}^{\infty} |mz+n|^{-\sigma} \leq 8\zeta(\sigma-1)r(z)^{-\sigma}\ (\sigma > 2)$$

Thus the right-hand side of (1) converges absolutely and uniformly on any compact subset of H. Specifically, $E_k(z; \chi, \psi)$ is a holomorphic function on H. Putting:

$$(3)\ \Gamma_0(L, M) = \left\{ \begin{bmatrix} a & b \\ c & d \end{bmatrix} \in SL_2(Z) \mid b \equiv 0 \bmod M, c \equiv 0 \bmod L \right\}$$

Then $\Gamma_0(L, M)$ is modular form. This is an important result that any graduate student working in analytic number theory should learn.

Now for Ramanujan's work on modular forms: Ramanujan studied the Fourier coefficients of meromorphic modular forms. In particular, Ramanujan studied the reciprocal of the weight 6 Eisenstein series. What it would be important for any number theory student to know is the important relationship between Ramanujan's work on partitions and his work on modular forms. We learned that Ramanujan and Hardy used the circle method to derive their asymptotic formula for the partition function p(n):

$$p(n) \sim \frac{exp(\pi\sqrt{2n/3})}{4n\sqrt{3}}$$

Rademacher later derived the exact formula:

$$p(n) = 2\pi(24n-1)^{-2/3} \sum_{k=1}^{\infty} \frac{A_k(n)}{k} \cdot I_{\frac{3}{2}}\left(\frac{\pi\sqrt{24n-1}}{6k}\right)$$

Where $I_l(x)$ denotes the I-Bessel function of order l and $A_k(n)$ denotes the Kloosterman sum.

What any student working in analytic number theory should know is that there is a rich connection between Ramanujan's work on partitions and his theory of modular forms. For example, Rademacher's exact formula comes from considering the partition function as a type of modular function. In particular, the function:

$$P(z) := \sum_{n=0}^{\infty} p(n) q^{n-\frac{1}{24}} = q^{-\frac{1}{24}} \prod_{n=1}^{\infty} \frac{1}{1-q^n}$$

is a weight -1/2 weakly modular form.

Finally, there are applications of Ramanujan's work to the field of physics. One area where Ramanujan's mathematics has applications is in the study of the moonshine phenomenon. There is an important relationship between work on moonshine and the Rogers-Ramanujan identities. The Rogers-Ramanujan identities are:

$$\sum_{n=0}^{\infty} \frac{q^{n^2}}{(q)_n} = \prod_{j=0}^{\infty} \frac{1}{(1-q^{5j+1})(1-q^{5j+4})},$$

$$\sum_{n=0}^{\infty} \frac{q^{n(n+1)}}{(q)_n} = \prod_{j=0}^{\infty} \frac{1}{(1-q^{5j+2})(1-q^{5j+3})}$$

where $(q)_n = (1-q)(1-q^2)\ldots(1-q^n)$ for n > 0 and $(q)_0 = 1$.

The Rogers-Ramanujan identities express two infinite product modular forms as number-theoretic q-series. The theory of Monstrous Moonshine asserts that the coefficients of the modular j-function are dimensions of virtual characters for the Monster group, the largest of the simple sporadic groups. The key to the relationship between the Rogers-Ramanujan identities and the theory of Monstrous Moonshine is, of course, modular forms.

## Conclusion

Ramanujan made many important contributions to infinite series, number theory, and complex analysis. Ramanjan's legacy will live on, going on to inspire generations of mathematicians to come.